\documentclass[12pt]{amsart}
\usepackage{amsfonts, amssymb}
\usepackage{amsmath}
\usepackage{amsthm}
\usepackage{cite}
\numberwithin{equation}{section}

\newtheorem{theorem}{Theorem}[section]

\newtheorem{corollary}[theorem]{Corollary}

\begin{document}

\title[$q$-difference equations of Malmquist type]{Zero order meromorphic solutions of $q$-difference equations of Malmquist type}

\author{Risto Korhonen}
\address{Department of Physics and Mathematics, University of Eastern Finland, P.O. Box 111, FI-80101 Joensuu, Finland,\newline
%
%
}
\email{risto.korhonen@uef.fi}

\author{Yueyang Zhang$^{*}$}
\footnote{*Corresponding author}
\address{School of Mathematics and Physics, University of Science and Technology Beijing, No.~30 Xueyuan Road, Haidian, Beijing, 100083, P.R. China}
\email{zhangyueyang@ustb.edu.cn}
\thanks{The second author is supported by the Fundamental Research Funds for the Central Universities~{(FRF-TP-19-055A1)} and a Project supported by the National Natural Science Foundation of China~{(12301091)}}



\date{\today}

\commby{}

\begin{abstract}
We consider the first order $q$-difference equation
\begin{equation}\tag{\dag}
f(qz)^n=R(z,f),
\end{equation}
where $q\not=0,1$ is a constant and $R(z,f)$ is rational in both arguments. When $|q|\not=1$, we show that, if $(\dag)$ has a zero order transcendental meromorphic solution, then $(\dag)$ reduces to a $q$-difference linear or Riccati equation, or to an equation that can be transformed to a $q$-difference Riccati equation. In the autonomous case, explicit meromorphic solutions of $(\dag)$ are presented. Given that $(\dag)$ can be transformed into a difference equation, we proceed to discuss the growth of the composite function $f(\omega(z))$, where $\omega(z)$ is an entire function satisfying $\omega(z+1)=q\omega(z)$, and demonstrate how the proposed difference Painlev\'e property, as discussed in the literature, applies for $q$-difference equations.
\end{abstract}

\maketitle


\section{Introduction}\label{intro} 

An ordinary differential equation is said to possess the Painlev\'{e} property if all of its solutions are single-valued about all movable singularities. It is widely believed that all ordinary differential equations that possess the Painlev\'{e} property are integrable. In~\cite{AblowitzHalburdHerbst2000}, Ablowitz, Halburd and Herbst introduced Nevanlinna theory (see, e.g., \cite{Hayman1964Meromorphic}) to study the integrability of difference equations and suggested that the existence of sufficiently many finite-order meromorphic solutions of a difference equation is a good candidate for a difference analogue of the Painlev\'e property. Halburd and Korhonen~\cite{halburdrk:LMS2006} implemented this idea to the second order difference equation
\begin{equation}\label{ahh_eq}
f(z+1)+f(z-1) = R(z,f),
\end{equation}
where $R(z,f)$ is rational in $f$ with meromorphic coefficients having slower growth than $f(z)$ in terms of Nevanlinna theory, and showed that if equation \eqref{ahh_eq} has a meromorphic solution $f$ of finite order, then either $f$ satisfies a difference Riccati equation, or a linear transformation of \eqref{ahh_eq} reduces it to one in a short list of difference equations which consists solely of difference Painlev\'e equations and equations related to them, linear equations and linearizable equations. It was shown that the finite-order condition of the proposed difference Painlev\'e property can be relaxed to hyper-order strictly less than one in \cite{halburdkt:14TAMS}, and later to minimal hyper-type (i.e. $\log T(r,f)=o(r)$ as $r\to\infty$) in~\cite{Korhonenkazhzh2020,zhengR:18}.

When detecting the integrability of $q$-difference equation, it is natural to consider the zero order meromorphic solutions. One of the main purposes of this paper is to display how the proposed difference Painlev\'e property works for $q$-difference equations. We shall consider the first order $q$-difference equation
\begin{equation}\label{first_order_q-ver1}
f(qz)^n=R(z,f),
\end{equation}
where $n\in\mathbb{N}$, $q\not=0,1$ is a constant and $R(z,f)$ is rational in both arguments. Equation \eqref{first_order_q-ver1} is the $q$-difference version of the first order difference equation
\begin{equation}\label{first_order_de_n}
g(z+1)^n=R(z,g),
\end{equation}
where $n\in \mathbb{N}$ and $R(z,g)$ is rational in $g$ with meromorphic coefficients having slower growth than $g(z)$. In fact, we may perform the transformation $z\to \omega(z)$ to equation \eqref{first_order_q-ver1} with an entire function $\omega(z)$ satisfying $\omega(z+1)=q\omega(z)$. We may choose $\omega(z)=\kappa(z)e^{z\log q}$ with an arbitrary entire function $\kappa(z)$ such that $\kappa(qz)=\kappa(z)$. Denoting $g(z)=f(\omega(z))$, then $g(z+1)=f(q\omega(z))$, and $g(z)$ is a meromorphic solution of the difference equation
\begin{equation}\label{first_order_q-ver5}
g(z+1)^n=\hat{R}(z,g),
\end{equation}
where $\hat{R}(z,g)$ is now rational in $g$ with coefficients that are rational in $\omega(z)$. According to a result of Edrei and Fuchs \cite{EdreiFuchs1953} (see also \cite[Lemma~2.6]{Hayman1964Meromorphic}), we have
\begin{equation}\label{first_order_q-ver6}
\begin{split}
\lim_{r\to\infty}\frac{T(r,f(\omega(z)))}{T(r,\omega(z))}=\infty.
\end{split}
\end{equation}
Since, by the Valiron--Mohon'ko identity \cite{Valiron1931,mohonko1971} (see also \cite{Laine1993}), the characteristic function of each coefficient of $\hat{R}(z,g)$, say $\hat{a}(z)$, satisfies $T(r,\hat{a}(z))=\hat{n}T(r,\omega(z))+O(\log r)$, where $\hat{n}$ is some integer, all coefficients of \eqref{first_order_q-ver5} have growth of type $o(T(r,g))$ when the solution $f(z)$ of \eqref{first_order_q-ver1} is transcendental. Under the condition that each coefficient of \eqref{first_order_de_n}, say $\tilde{a}(z)$, satisfies $T(r,\tilde{a}(z+k))=o(T(r,g))$ for any finite positive integer $k$, where $r\to\infty$ outside an exceptional set of finite logarithmic measure (i.e. $\int_{E}dr/r<\infty$), the present authors \cite{Korhonenzhang2020,zhangkorhonen2022,Korhonenzhang2023} used Yamanoi's Second Main Theorem \cite{yamanoi:04,yamanoi:05} to provide a classification for equation \eqref{first_order_de_n}, including a list of $33$ equations of canonical form. In particular, in the autonomous case all these equations can be solved in terms of elliptic functions, exponential type functions or functions which are solutions to a certain autonomous first-order difference equation. These results provide a complete difference analogue of the results on Malmquist type differential equations $f'(z)^n=R(z,f)$ due to Steinmetz \cite{steinmetz:78} and Bank and Kaufman \cite{bank1980growth}; see also \cite[Chapter~10]{Laine1993}.

Now, if equation \eqref{first_order_q-ver1} has a zero order transcendental meromorphic solution $f$, then it follows from an estimate for the proximity function $m(r,f(qz)/f(z))$ in \cite{BarnettHalburd} (see also \cite{zhangjkorhonen2010}) and the Valiron--Mohon'ko identity that
\begin{equation*}
\begin{split}
\deg_f(R(z,f))T(r,f) &= T(r,R(z,f)) + O(\log r)\\
&= T(r,f(qz)^n) + O(\log r) \\
&= nT(r,f(qz)) + O(\log r) \\
&= nT(r,f) + o(T(r,f)),
\end{split}
\end{equation*}
where $r\to\infty$ on a set of logarithmic density~1 (i.e. $\lim_{r\to\infty}(\int_{E\cap [0,r]}dt/t)/\log r=1$). This implies that $\deg_f(R(z,f))=n$. Together with previous discussions, we now discard the assumption that $f(z)$ is of zero order and conclude from the results in \cite{Korhonenzhang2020,Korhonenzhang2023} for equation \eqref{first_order_de_n} with $\deg_f(R(z,f))=n$ that: If the $q$-difference equation \eqref{first_order_q-ver1} has a transcendental meromorphic solution $f$, then either $f$ satisfies a $q$-difference linear or Riccati equation
\begin{eqnarray}
f(qz)&=& a_1(z)f(z)+a_2(z),\label{q-lineareq}\\
f(qz)&=&\frac{b_1(z)f(z)+b_2(z)}{f(z)+b_3(z)}\label{q-driccati},
\end{eqnarray}
where $a_i(z)$ and $b_j(z)$ are rational functions, or, by implementing a transformation $f\rightarrow \alpha f$ or $f\rightarrow 1/(\alpha f)$ with an algebraic function $\alpha$ of degree at most $3$, \eqref{first_order_q-ver1} reduces to one of the following equations
\begin{eqnarray}
f(qz)^2 &=& 1-f(z)^2,\label{yanagiharaeq11 co}\\
f(qz)^2 &=& 1-\left(\frac{\delta(z) f(z)-1}{f(z)-\delta(z)}\right)^2,\label{yanagiharaeq12 co}\\
f(qz)^2 &=& 1-\left(\frac{f(z)+3}{f(z)-1}\right)^{2},\label{yanagiharaeq13 co}\\
f(qz)^2 &=& \frac{f(z)^2-\kappa(z)}{f(z)^2-1},\label{yanagiharaeq14 co}\\
f(qz)^3 &=& 1-f(z)^{-3},\label{yanagiharaeq15 co}
\end{eqnarray}
where $\delta(z)\not\equiv\pm1$ is an algebraic function of degree at most $2$ and $\kappa(z)\not\equiv0,1$ is a rational function such that $\kappa(qz)=\kappa(z)$, or to one of the following equations
\begin{eqnarray}
f(qz)^2 &=& \delta_1(z)(f(z)^2-1), \label{intro_eq_list2_1}\\
f(qz)^2 &=& \delta_2(z)(1-f(z)^{-2}), \label{intro_eq_list2_2}\\
f(qz)^2 &=& \frac{\delta_3(qz)f(z)^2-1}{f(z)^2-1}, \label{intro_eq_list2_3}\\
f(qz)^2 &=& \theta\frac{f(z)^2-\delta_4(z)f(z)+1}{f(z)^2+\delta_4(z)f(z)+1}, \label{intro_eq_list2_4}\\
f(qz)^3 &=& 1-f(z)^3,\label{intro_eq_list2_5}\\
f(qz)^2 &=& \frac{1}{2}\frac{(1+\delta_5(z))^2}{1+\delta_5(z)^2}\frac{(f(z)-1)(f(z)-\delta_5(z)^2)}{(f(z)-\delta_5(z))^2},  \label{intro_eq_list2_6}
\end{eqnarray}
where $\theta=\pm1$ and $\delta_1$, $\delta_2$, $\delta_3$, $\delta_4$ and $\delta_5$ are rational or algebraic functions, each of which satisfying
\begin{eqnarray}
\delta_1(qz)(\delta_1(z)+1)+1 &=& 0, \label{intro_eq_list2_1 co}\\
\delta_2(qz)\delta_2(z) &=& \delta_2(qz)+\delta_2(z), \label{intro_eq_list2_2 co}\\
\delta_3(qz)\delta_3(z) &=& 1, \label{intro_eq_list2_3 co}\\
\delta_4(qz)(\delta_4(z)-4) &=& 2(1-\theta)\delta_4(z)-8(1+\theta), \label{intro_eq_list2_4 co}\\
8\delta_5^4(qz)(\delta_5^2(z)+1)\delta_5(z) &=& (\delta_5(z)+1)^4,\label{intro_eq_list2_5 co}
\end{eqnarray}
respectively. From the proof in \cite{Korhonenzhang2020,Korhonenzhang2023} we know that solutions of \eqref{yanagiharaeq13 co} are still meromorphic after the transformation $f\to \alpha f$ or $f\to 1/(\alpha f)$. Concerning zero order meromorphic solutions of \eqref{first_order_q-ver1}, we shall prove the following

\begin{theorem}\label{maintheorem-1-1}
Suppose that $|q|\not=1$ in \eqref{first_order_q-ver1}. If equation \eqref{first_order_q-ver1} has a zero order transcendental meromorphic solution, then either $f$ satisfies \eqref{q-lineareq} or \eqref{q-driccati} with rational coefficients, or, by a transformation $f\rightarrow \alpha f$ or $f\rightarrow 1/(\alpha f)$ with an algebraic function $\alpha$ of degree at most $2$, \eqref{first_order_q-ver1} reduces to \eqref{yanagiharaeq12 co}.
\end{theorem}

By combining the results in the first two parts of Section~\ref{Derivations and existence}, we have actually shown that, if $|q|\not=1$, then equation \eqref{first_order_q-ver1} can have transcendental meromorphic solutions only in the case of the four equations \eqref{q-lineareq}, \eqref{q-driccati}, \eqref{yanagiharaeq12 co} and \eqref{intro_eq_list2_6}.

The rather simple $q$-difference equation $f(z)-a(z)f(qz)=0$, where $a(z)$ is a nonconstant polynomial, has a zero order transcendental entire solution~\cite[Theorem~1]{Bergweilerhayman2004}. In fact, Bergweiler, Ishizaki and Yanagihara~\cite{BergweilerIshizaki1998} have previously shown that all meromorphic solutions of the $q$-difference equation $\sum_{j=0}^ma_j(z)g(q^iz)=Q(z)$, where $q\in \mathbb{\mathbb{C}}$, $0<|q|<1$, and $Q$ and the $a_j$ are polynomials, are of zero order. They also showed that, under a suitable assumption on $a_0(z)$ and $a_1(z)$, the $q$-difference equation $g(q^2z)+a_1(z)g(qz)+a_0(z)g(z)=0$ has a transcendental meromorphic solution $g(z)$. Letting $f(z)=g(qz)/g(z)$, then $f(z)$ satisfies the $q$-difference Riccati equation $f(qz)=-[a_1(z)f(z)+a_0(z)]/f(z)$. With this result we may construct zero order meromorphic solutions of equation \eqref{yanagiharaeq12 co} for nonconstant $\delta$. From \cite{Korhonenzhang2020} we know that the solution $f(z)$ of equation \eqref{yanagiharaeq12 co} is represented as $f(z)=(\gamma(z)+\gamma(z)^{-1})/2$, where $\gamma(z)$ satisfies a $q$-difference Riccati equation
\begin{equation}\label{riccati  conlim11fua1fat0000}
\gamma(qz)=\left\{-\theta\frac{(\pm i\delta(z)-\sqrt{1-\delta(z)^2})\gamma(z)\pm i}{\gamma(z)-\delta(z)\pm i\sqrt{1-\delta(z)^2}}\right\}^{\theta}, \quad \theta=\pm1.
\end{equation}
Consider the case $\theta=-1$ and $\pm$ is chosen to be $-$. We choose $\delta=\frac{1}{2}(\alpha(z)+\alpha(z)^{-1})$ with a suitable polynomial $\alpha(z)$ so that $a_1(z)=-(\alpha(qz)+i\alpha(z))$ and $a_0(z)=2i\alpha(z)$. Then, by letting $\Gamma(z)=\gamma(z)+\alpha(z)$, we have
\begin{equation}\label{riccati  conlim11fua1fat0000ju}
\Gamma(qz)=\frac{(\alpha(qz)+i\alpha(z))\Gamma(z)-2i\alpha(z)^{2}}{\Gamma(z)}.
\end{equation}
In the autonomous case, the meromorphic solutions of the $q$-difference linear and Riccati equations \eqref{q-lineareq} and \eqref{q-driccati} have been clearly characterised in~\cite{Goldstein1970}. In Section~\ref{Derivations and existence} we shall present some explicit meromorphic solutions for these two equations.

In the proof of Theorem~\ref{maintheorem-1-1}, we show that $\delta_5(z)$ in \eqref{intro_eq_list2_6} must be a constant. Then we can determine all nonconstant meromorphic solutions of equation \eqref{intro_eq_list2_6}, which are actually of positive even order. However, in the case of the seven equations \eqref{yanagiharaeq14 co}--\eqref{intro_eq_list2_5}, the existence of transcendental meromorphic solutions does not necessarily imply that the coefficients of \eqref{first_order_q-ver1} are constants. Nonetheless, when all coefficients of equation \eqref{first_order_q-ver1} are constants, we may remove the condition $|q|\not=1$ in Theorem~\ref{maintheorem-1-1} and obtain the following

\begin{corollary}\label{corollary-1-2}
If equation \eqref{first_order_q-ver1} with constant coefficients has a nonconstant zero order meromorphic solution, then either $f$ satisfies \eqref{q-lineareq} or \eqref{q-driccati} with constant coefficients, or, by a transformation $f\rightarrow \alpha f$ or $f\rightarrow 1/(\alpha f)$ with a constant $\alpha$, \eqref{first_order_q-ver1} reduces to equations \eqref{yanagiharaeq11 co}, \eqref{yanagiharaeq12 co}, \eqref{yanagiharaeq13 co} with constant coefficients.
\end{corollary}

We remark that, when $|q|=1$, equation \eqref{first_order_q-ver1} with nonconstant rational coefficients can also have zero order transcendental meromorphic solutions in the case of equation \eqref{yanagiharaeq11 co}.
From \cite{Korhonenzhang2020} we know that the solution $f(z)$ of \eqref{yanagiharaeq11 co} is represented as $f(z)=(h(z)+h(z)^{-1})/2$, where $h(z)$ satisfies $h(qz)=ih(z)^{\theta}$, and $\theta=\pm1$. If $f(z)$ is written as $f(z)=g(z)/\alpha(z)$ for an algebraic function $\alpha(z)$ such that $\alpha(z)^2$ is a rational function and $\alpha(qz)^2=\alpha(z)^2$, then we see that $h(z)$ takes the form $h(z)=w(z)/\alpha(z)$ for a meromorphic function $w(z)$ such that $w(qz)/\alpha(z)=i(w(z)/\alpha(z))^{\theta}$. This equation can have zero order transcendental meromorphic solutions when $|q|=1$.

This paper is structured as follows. In Section~\ref{Derivations and existence}, we first present meromorphic solutions for equation \eqref{first_order_q-ver1} in the autonomous case. The $11$ equations \eqref{yanagiharaeq11 co}--\eqref{intro_eq_list2_6} are solved in terms of rational functions or elliptic functions. Together with these solutions, we then provide a proof for Theorem~\ref{maintheorem-1-1}. Since equation \eqref{first_order_q-ver1} can be transformed to an equation of the form in \eqref{first_order_de_n}, in the final part of Section~\ref{Derivations and existence}, we discuss the relationship between solutions of equation \eqref{first_order_q-ver1} and \eqref{first_order_de_n} by analysing the growth of the composite function $g(z)=f(\omega(z))$, where $\omega(z)$ is an entire function such that $\omega(z+1)=q\omega(z)$. We point out that the existence of sufficiently many zero order meromorphic solutions of the $q$-difference equation \eqref{first_order_q-ver1} implies the existence of sufficiently many minimal hyper-type meromorphic solutions of the difference equation \eqref{first_order_de_n}. Finally, in Section~\ref{Concluding remarks} we make two comments on the integrability of $q$-difference equations.

\section{Zero order meromorphic solutions of equation \eqref{first_order_q-ver1}}\label{Derivations and existence}

\subsection{Meromorphic solutions of \eqref{first_order_q-ver1} in the autonomous case}\label{examples} 

In this section we present meromorphic solutions of equation \eqref{first_order_q-ver1} in the autonomous case. In particular, we will show that the five equations \eqref{q-lineareq}--\eqref{yanagiharaeq13 co} indeed have zero order meromorphic solutions and nonconstant meromorphic solutions of the eight equations \eqref{yanagiharaeq14 co}--\eqref{intro_eq_list2_6} have even positive order or infinite order.

For the autonomous $q$-difference linear equation \eqref{q-lineareq}: $f(qz)=af(z)+b$, where $a,b$ are constants, the general solution is $f(z)=z^{\log a/\log q}+b/(1-a)$ when $|a|\not=1$; or $f(z)=\kappa(z)z^{\log a/\log q}+b/(1-a)$ when $|a|=1$ but $a\not=1$, where $\kappa(z)$ is an arbitrary function such that $\kappa(qz)=\kappa(z)$; or $f(z)=\kappa(z)+b\log z/\log q$ when $a=1$, where $\kappa(z)$ is an arbitrary function such that $\kappa(qz)=\kappa(z)$. When $a\not=1$, if the solution $f(z)$ is meromorphic, then the constant $q$ should satisfy $\log a/\log q=k$ for an integer $k$; when $a=1$, if the solution $f(z)$ is meromorphic, then the coefficient $b$ should vanish. In particular, we may choose $\kappa(z)$ to be a zero order transcendental meromorphic (entire) function. We remark that, if $|q|\not=1$ or $|q|=1$ and $q=e^{2\pi i\theta}$ for an irrational number $\theta\in[0,1)$, then the equation $\kappa(qz)=\kappa(z)$ does not have any non-rational zero order meromorphic solutions. This is due to the fact that any possible zeros or poles of the solution accumulate (so they must be Picard exceptional) and a nonconstant zero order meromorphic function has at most one Picard exceptional value.

For the autonomous $q$-difference Riccati equation \eqref{q-driccati}, by a result of Ishizaki~\cite{Ishizaki2017}, after doing certain transformations we obtain the $q$-difference Riccati equation
\begin{equation}\label{q-riccati  1}
\begin{split}
f(qz)=\frac{f(z)+A}{1-f(z)}
\end{split}
\end{equation}
or
\begin{equation}\label{q-riccati  2}
\begin{split}
f(qz)=\frac{B}{f(z)},
\end{split}
\end{equation}
where $A$ and $B\not=0$ are two constants. For equation \eqref{q-riccati  1}, let $a_1$ and $a_2$ be two nonzero constants and $b=-2a_1a_2/(a_1^2+a_2^2)$. When $A=-b^2$ and $a$ is a constant such that $e^{ia}=(a_1-a_2)^2/(a_1+a_2)^2$, we have the solution
\begin{equation}\label{riccati  5}
f(z)=-\left(\frac{2a_1a_2}{a_1^2+a_2^2}\right)\frac{\kappa(z)z^{a/\log q}-a_1}{\kappa(z)z^{a/\log q}+a_2},
\end{equation}
where $\kappa(z)$ is an arbitrary function such that $\kappa(qz)=\kappa(z)$. For equation \eqref{q-riccati  2}, letting $a_3$ be a nonzero constant, then we have the solution $f(z)=a_3Bw(z)/w(qz)$, where $w(z)$ is an arbitrary entire function such that $w(q^2z)=a_3^2B^2w(z)$. Obviously, we may choose $\kappa(z)$ or $w(z)$ to be a transcendental meromorphic function of zero order, as mentioned in previous discussions on equation \eqref{q-lineareq}.

Look at the three equations \eqref{yanagiharaeq11 co}, \eqref{yanagiharaeq12 co} and \eqref{yanagiharaeq13 co}. They are so-called Fermat $q$-difference equations. Recall that a Fermat equation is a function analogue of the Fermat diophantine equation $x^n+y^n=1$, i.e., $x(z)^n+y(z)^n=1$, where $n\geq2$ is an integer. Meromorphic solutions to Fermat equations have been clearly characterized; see \cite{Baker1966,Gross1966erratum,Gross1966}, for example. Solutions of the three equations \eqref{yanagiharaeq11 co}, \eqref{yanagiharaeq12 co} and \eqref{yanagiharaeq13 co} are solved in a similar way as in \cite{Korhonenzhang2020}. For equation \eqref{yanagiharaeq11 co}, the solution $f(z)$ is represented as $f(z)=(\beta(z)+\beta(z)^{-1})/2$, where $\beta(z)$ satisfies $\beta(qz)=i\beta(z)^{\pm1}$. For equation \eqref{yanagiharaeq12 co}, the solution $f(z)$ is represented as $f(z)=(\gamma(z)+\gamma(z)^{-1})/2$, where $\gamma(z)$ satisfies
    \begin{equation}\label{riccati  spe1}
    \gamma(qz)=\left\{-\theta\frac{(i\delta-\sqrt{1-\delta^2})\gamma(z)+i}{\gamma(z)-\delta+i\sqrt{1-\delta^2}}\right\}^{\theta}, \quad \theta=\pm1,
    \end{equation}
where $\delta$ is a constant. Moreover, based on the discussions in \cite{Korhonenzhang2023} we know that equation \eqref{riccati  spe1} can be transformed into equation \eqref{q-riccati  1} when $2\delta^2\not=1$ and can be transformed into equation \eqref{q-riccati  2} when $2\delta^2=1$. For equation \eqref{yanagiharaeq13 co}, the solution $f(z)$ is represented as $f(z)=\frac{8\lambda(z)^2-(\lambda(z)^2+1)^2}{(\lambda(z)^2+1)^2}$, where $\lambda(z)$ satisfies
    \begin{equation}\label{riccati  spe2}
    \lambda(qz)=\left\{-\theta\frac{-(1+\sqrt{2})\lambda(z)+i}{\lambda(z)-i+i\sqrt{2}}\right\}^{\theta}, \quad \theta=\pm1.
    \end{equation}
Equation \eqref{riccati  spe2} can be transformed into equation \eqref{q-riccati  1}. Thus meromorphic solutions of the three equations \eqref{yanagiharaeq11 co}, \eqref{yanagiharaeq12 co} and \eqref{yanagiharaeq13 co} are solved in terms of functions which are solutions of certain $q$-difference linear or Riccati equation.

Now look at the six equations \eqref{yanagiharaeq14 co}, \eqref{intro_eq_list2_1}, \eqref{intro_eq_list2_2}, \eqref{intro_eq_list2_3}, \eqref{intro_eq_list2_4} and \eqref{intro_eq_list2_6}. The nonconstant solutions of these equation are Jacobi elliptic functions composed with entire functions.

We first take the equation \eqref{intro_eq_list2_6} as an example and show that all nonconstant meromorphic solutions of this equation have even positive order. From the proof in \cite{Korhonenzhang2023} we know that a meromorphic solution $f(z)$ of equation \eqref{intro_eq_list2_6} is twofold ramified over each of $\pm1,\pm \delta_5^2$. Then there exists an entire function $\phi(z)$ such that $f(z)$ is written as $f(z)=\text{sn}(\phi(z))$, where $\text{sn}(\phi)=\text{sn}(\phi,1/\delta_5^2)$ is the Jacobi elliptic function with modulus $1/\delta_5^2$ and satisfies the first order differential equation $\text{sn}'(\phi)^2=(1-\text{sn}(\phi)^2)(1-\text{sn}(\phi)^2/\delta_5^4)$. Let $z_0$ be a point such that $f(z_0)=\text{sn}(\phi(z_0))=1$. It follows from \eqref{intro_eq_list2_6} that $f(qz_0)=\text{sn}(\phi(qz_0))=0$. Using the differential equation $\text{sn}'(\phi)^2=(1-\text{sn}(\phi)^2)(1-\text{sn}(\phi)^2/\delta_5^4)$ to compute the Taylor series for $\text{sn}(\phi)$ and $\text{sn}(\phi(qz))$ around the point $z_0$, respectively, we get
    \begin{equation}\label{Poiu1}
    \begin{split}
    \text{sn}(\phi(z))&=1-\frac{1}{2}\frac{\delta_5^4-1}{\delta_5^4}(\phi(z)-\phi(z_0))^2+\cdots\\
    &=1-\frac{1}{2}\frac{\delta_5^4-1}{\delta_5^4}\phi'(z_0)^2(z-z_0)^2+\cdots
    \end{split}
    \end{equation}
and
    \begin{equation}\label{Poiu2}
    \begin{split}
   \text{sn}(\phi(qz))=\phi(qz)-\phi(qz_0)+\cdots=q\phi'(qz_0)(z-z_0)+\cdots.
    \end{split}
    \end{equation}
By substituting the above two expressions into \eqref{intro_eq_list2_6} and then comparing the second-degree terms on both sides of the resulting equation, we find
    \begin{equation}\label{Poiu3}
    \begin{split}
    q^2\phi'(qz_0)^2 =\frac{1}{4}\frac{(1+\delta_5)^4}{\delta_5^4}\phi'(z_0)^2.
    \end{split}
    \end{equation}
Define $G(z):=q^2\phi'(qz)^2-\frac{1}{4}\frac{(1+\delta_5)^4}{\delta_5^4}\phi'(z)^2$. Similarly as in \cite{Korhonenzhang2023}, we can use \eqref{first_order_q-ver6} to prove that equation \eqref{Poiu3} holds for all $z\in \mathbb{C}$, so that $G(z)\equiv0$. It follows by integration that
    \begin{equation}\label{Poiu4}
    \begin{split}
    \phi(qz)=b_1\phi(z)+b_2, \ \ \ b_1=\pm \frac{1}{2}\frac{(1+\delta_5)^2}{\delta_5^2}.
    \end{split}
    \end{equation}
where $b_2$ is a constant. (The coefficient $b_1$ is different from the one in equation (2.34) in \cite{Korhonenzhang2023}. The calculation there contains an error.) Since $\delta_5\not=\pm i$, we see that $|b_1|\not=1$ and thus $\phi$ is a polynomial of the form $B_1z^k+B_2$ with $k$ being a positive integer such that $q^k=b_1$. Since $\text{sn}(z)$ has order~2, then by a result of Bergweiler \cite[Satz~3.2]{BergweilerJankVol1984} we have that the order of $f(z)=\text{sn}(\phi(z))$ is equal to~$2k$.

Similarly, for each of the five equations \eqref{yanagiharaeq14 co}, \eqref{intro_eq_list2_1}, \eqref{intro_eq_list2_2}, \eqref{intro_eq_list2_3} and \eqref{intro_eq_list2_4}, we can show that the nonconstant solutions are written as $f(z)=\text{sn}(\phi_1)=\text{sn}(\phi_1,\delta)$, where $\text{sn}(\phi_1)$ is the Jacobi elliptic function with a certain modulus $\delta$ and an entire function $\phi_1$ satisfying a $q$-difference linear equation
\begin{equation}\label{intro_eq_list2_6 fuqjh[kokut}
\phi_1(qz) =b_3\phi_1(z)+b_4,
\end{equation}
where $b_3$ and $b_4$ are two constants such that $|b_3|=1$. In particular, for equation \eqref{yanagiharaeq14 co}, we have $b_3=\pm 1$. We show that $b_3=1$ is impossible. Otherwise, from the previous discussions on equation \eqref{q-lineareq} we see that $b_4=0$. Recall the Maclaurin series for $\text{sn}\ \varepsilon$:
\begin{equation}\label{MYS2 fu5jh5}
\begin{split}
\text{sn}\ \varepsilon&=\varepsilon-(1+\delta_0^2)\frac{\varepsilon^3}{3!}+(1+14\delta_0^2+\delta_0^4)\frac{\varepsilon^5}{5!}+\cdots.
\end{split}
\end{equation}
Since $\kappa\not=0,1$ and $\text{sn}(\phi_1(qz))=\text{sn}(\phi_1(z))$, then substitution into equation \eqref{yanagiharaeq14 co} yields a contradiction. If $\phi_1$ is a polynomial, then $f(z)$ has finite order $2k_1$ for some positive integer $k_1$; if $\phi_1$ is transcendental, then $f(z)$ has infinite order by a result of \cite[Theorem~2.9]{Hayman1964Meromorphic}.

Then look at the third degree Fermat $q$-difference equations \eqref{yanagiharaeq15 co} and \eqref{intro_eq_list2_5}. From \cite{Baker1966,Gross1966erratum,Gross1966} we know that all meromorphic solutions of the equation $x(z)^3+y(z)^3=1$ can be represented as: $x=H(\varphi)$, $y=\eta G(\varphi)=\eta H(-\varphi)=H(-\eta^2\varphi)$, where $\varphi=\varphi(z)$ is an entire function and $\eta$ is a cubic root of~1, and
\begin{equation}\label{Weierstrassellptic0}
H(z)=\frac{1+\wp'(z)/\sqrt{3}}{2\wp(z)}, \quad G(z)=\frac{1-\wp'(z)/\sqrt{3}}{2\wp(z)},
\end{equation}
where $\wp(z)$ is the particular Weierstrass elliptic function satisfying $\wp'(z)^2=4\wp(z)^3-1$. For equation \eqref{yanagiharaeq15 co}, we have $f(qz)=H(\varphi_1)$ and $f(z)^{-1}=\eta G(\varphi_1)$ with an entire function $\varphi_1(z)$. An elementary series analysis on $f(z)$ as for equation \eqref{intro_eq_list2_6} shows that $\varphi_1(z)$ satisfies a first order $q$-difference linear equation
\begin{equation}\label{intro_eq_list2_6 fuqjh[kokutqw}
\varphi_1(qz)=c_1\varphi_1(z)+c_2,
\end{equation}
where $c_1$ and $c_2$ are constants such that $c_1^3=1$. Moreover, by the same arguments as for equation \eqref{yanagiharaeq14 co}, we may show that $c_1\not=1$. For equation \eqref{intro_eq_list2_5}, we have $f(qz)=H(\varphi_2)$ and $f(z)=\eta G(\varphi_2)$ with an entire function $\varphi_2$. An elementary series analysis on $f(z)$ as for equation \eqref{intro_eq_list2_6} shows that $\varphi_2(z)$ satisfies a first order $q$-difference linear equation
\begin{equation}\label{intro_eq_list2_6 fuqjh[kokutqwg}
\varphi_2(qz)=c_3\varphi_2(z)+c_4,
\end{equation}
where $c_3$ and $c_4$ are constants such that $c_2^3=-1$. Note that the Weierstrass elliptic function $\wp(z)$ has order $2$. If $\varphi_1$ (or $\varphi_2$) is a polynomial, then the solution $f(z)$ of \eqref{yanagiharaeq15 co} (or \eqref{intro_eq_list2_5}) has finite order $2m$ for some positive integer $m$; if $\varphi_1$ (or $\varphi_2$) is transcendental, then the solution $f(z)$ of \eqref{yanagiharaeq15 co} (or \eqref{intro_eq_list2_5}) has infinite order.

Finally, we show how to determine the nonconstant rational solutions of the five equations \eqref{q-lineareq}--\eqref{yanagiharaeq13 co} in the autonomous case. From the remarks in \cite[Section~2]{GundersenHeittokangas2002} we know that the autonomous $q$-difference linear and Riccati equations \eqref{q-lineareq} and \eqref{q-driccati} have only rational solutions when $|q|\not=1$. We see that the autonomous $q$-difference linear equation \eqref{q-lineareq} can have nonconstant rational solutions in both of the two cases $|q|\not=1$ and $|q|=1$. Below we shall show that equation \eqref{q-driccati} has no nonconstant rational solutions when $|q|\not=1$.

We only need to consider the two $q$-difference equations \eqref{q-riccati  1} and \eqref{q-riccati  2}. Suppose that equation \eqref{q-riccati  1} has a nonconstant rational solution $f(z)$. We may write $f(z)$ as $f(z)=P(z)/Q(z)$, where $P(z)$ and $Q(z)$ are two coprime polynomials of degrees $s_1$ and $s_2$ respectively. Denote the leading coefficients of $P(z)$ and $Q(z)$ by $a$ and $1$, respectively. Note that $A\not=-1$. When $A\not=0$, it is easy to see that $P(z)$ and $Q(z)$ have the same degrees, say $s$, and also that $P(0)\not=0$ and $Q(0)\not=0$. Since $a\not=1$, then by substituting $f(z)=P(z)/Q(z)$ into \eqref{q-riccati  1}, we may have
\begin{equation}\label{q-riccati  3}
\begin{split}
q^{-s}P(qz)&=(1-a)^{-1}[P(z)+AQ(z)],\\
q^{-s}Q(qz)&=(1-a)^{-1}[-P(z)+Q(z)].
\end{split}
\end{equation}
From the above two equations we get
\begin{equation}\label{q-riccati  3 g}
\begin{split}
P(q^2z)=2(1-a)^{-1}q^sP(qz)-(1-a)^{-2}q^{2s}(1+A)P(z).
\end{split}
\end{equation}
We write the formal form of $P(z)$:
\begin{equation}\label{q-riccati  3 ggy}
\begin{split}
P(z)=a_{s}z^{s}+a_{s-1}z^{s-1}+\cdots+a_0,
\end{split}
\end{equation}
where $a_{s}=a$ and $a_0\not=0$. By substituting $P(z)$ into \eqref{q-riccati  3 g} and then comparing the coefficients of the resulting polynomials on both sides, we get
\begin{equation}\label{q-riccati  3 ggy1}
\begin{split}
a_{s-i}=[2(1-a)^{-1}q^i-(1-a)^{-2}(1+A)q^{2i}]a_{s-i}, \ \ i=0,1,\cdots,s.
\end{split}
\end{equation}
By the two equations in \eqref{q-riccati  3 ggy1} in the case $i=0$ and $i=s$, we obtain $q^s=1$. We see that $a_s=a$ is determined by the equation $2(1-a)^{-1}-(1-a)^{-2}(1+A)=1$. Further, denoting by $\hat{s}$ the largest factor of $s$ such that $q^{\hat{s}}=1$, say $s/\hat{s}=\mathfrak{m}$, then the coefficients $a_{\hat{s}},a_{2\hat{s}},\cdots,a_{\mathfrak{m}\hat{s}}$ can be chosen arbitrarily while the other coefficients vanish. Then we can also give the form of $Q(z)$ using the identities in \eqref{q-riccati  3}. On the other hand, when $A=0$, it is easy to see that $s_1\leq s_2$. If $s_1=s_2$, then we have equations in \eqref{q-riccati  3} with $A=0$ and by comparing the leading coefficients of the polynomials on both sides of the first equation in \eqref{q-riccati  3} we get $a=0$, a contradiction. If $s_1<s_2$, then we have the following two equations
\begin{equation}\label{q-riccati  4}
\begin{split}
q^{-s}P(qz)&=P(z),\\
q^{-s}Q(qz)&=-P(z)+Q(z),
\end{split}
\end{equation}
from which we get $P(qz)=q^{s}P(z)$ and $Q(q^2z)=q^{s}Q(qz)-(q^{2s}-q^{s})Q(z)$. Since $Q(0)\not=0$, then by the same arguments as in the case $A\not=0$ we get $q^s=1$. However, it follows that $Q(qz)=Q(z)$, which is impossible.

Similarly, for equation \eqref{q-riccati  2}, substitution of $f(z)=P(z)/Q(z)$, where $P(z)$ and $Q(z)$ are two coprime polynomials of the same degree $s$ and that $P(0)\not=0$ and $Q(0)\not=0$, gives
\begin{equation}\label{q-riccati  5}
\begin{split}
q^{-s}P(qz)&=a^{-1}BQ(z),\\
q^{-s}Q(qz)&=a^{-1}P(z).
\end{split}
\end{equation}
It follows that $P(q^2z)=a^{-2}q^{2s}BP(z)$, which yields that $q^{2s}=1$. Similarly as for equation \eqref{q-riccati  1}, explicit forms for $P(z)$ and $Q(z)$ can be obtained from the system of two equations
in \eqref{q-riccati  5}. Then all nonconstant rational solutions of the five equations \eqref{yanagiharaeq11 co}--\eqref{yanagiharaeq13 co} with constant coefficients can be obtained after doing some transformations.

\subsection{Proof of Theorem~\ref{maintheorem-1-1}}\label{An elementary derivation} 

From the discussions in the introduction, we only need to consider the $13$ equations \eqref{q-lineareq}--\eqref{intro_eq_list2_6}. Below we first consider the four equations \eqref{intro_eq_list2_1}, \eqref{intro_eq_list2_2}, \eqref{intro_eq_list2_3} and \eqref{intro_eq_list2_4}.

From the discussions in \cite{Korhonenzhang2020} we know that, if we let $w=f+1/f$ for a solution $f$ of equation \eqref{intro_eq_list2_4}, then $w$ satisfies either an equation of the form in \eqref{intro_eq_list2_1} when $\theta=1$ or an equation of the form in \eqref{intro_eq_list2_2} when $\theta=-1$. Moreover, by the Valiron--Mohon'ko identity the characteristic function of $w$ satisfies $T(r,w)=2T(r,f)+O(1)$. Further, for equation \eqref{intro_eq_list2_2}, if we define a function $g(z)$ to be such that $-g(z)^2=f^2(z)-1$, then $g(z)$ satisfies
\begin{equation}\label{symm 4}
\begin{split}
g(qz)^2=[\hat{\delta}_2(qz)g(z)^2-1]/(g(z)^2-1),
\end{split}
\end{equation}
where $\hat{\delta}_2(z)$ satisfies $\hat{\delta}_2(qz)\hat{\delta}_2(z)=1$. This is just equation \eqref{intro_eq_list2_3}. Moreover, by the Valiron--Mohon'ko identity the characteristic function of $f(z)^2$ satisfies $T(r,f^2)=T(r,g^2)+O(1)$. Therefore, we only need to consider the two equations \eqref{intro_eq_list2_1} and \eqref{intro_eq_list2_2}.

Consider equation \eqref{intro_eq_list2_1}. In general, under the assumptions of Theorem~\ref{maintheorem-1-1} the solution $f(z)$ may have some algebraically branched points, but $f(z)^2$ is a meromorphic function. By iterating equation \eqref{intro_eq_list2_1} together with the relation in \eqref{intro_eq_list2_1 co} we find that
    \begin{equation}\label{iteration1}
    \begin{split}
    f(q^2z)^2&=-\frac{\delta_1}{1+\delta_1}f^2(z)+1,\\
    f(q^3z)^2&=f(z)^2.
    \end{split}
    \end{equation}
Then, at the origin $f(z)^2$ has a series expansion
    \begin{equation}\label{expansion1}
    \begin{split}
    f(z)^2=a_{-k}z^{-k}+\cdots+a_0+a_1z+\cdots.
    \end{split}
    \end{equation}
It follows that
    \begin{equation}\label{expansion2}
    \begin{split}
    f(q^3z)^2=a_{-k}q^{-3k}z^{-k}+\cdots+a_0+a_1q^3z+\cdots
    \end{split}
    \end{equation}
at the origin. From the above two series we conclude that $q^{3k}=1$ for some integer $k$, which is a contradiction to our assumption. Also, for equation \eqref{intro_eq_list2_2}, by iteration together with the relation in \eqref{intro_eq_list2_2 co} we easily find that $f(q^4z)^2=f(z)^2$. Then by the same arguments as for equation \eqref{intro_eq_list2_1} we have a contradiction. From the above reasoning, we conclude that in the case of equations \eqref{intro_eq_list2_1}, \eqref{intro_eq_list2_2}, \eqref{intro_eq_list2_3} and \eqref{intro_eq_list2_4}, equation \eqref{first_order_q-ver1} cannot have any zero order transcendental meromorphic solutions when $|q|\not=1$. In fact, when $|q|\not=1$, we see that these equations cannot have any transcendental meromorphic solutions.

By the the same arguments as for equation \eqref{intro_eq_list2_1} we can show that, in the case of the four equations \eqref{yanagiharaeq11 co}, \eqref{yanagiharaeq14 co}, \eqref{yanagiharaeq15 co} and \eqref{intro_eq_list2_5}, equation \eqref{first_order_q-ver1} cannot have any transcendental meromorphic solutions when $|q|\not=1$.

For equation \eqref{intro_eq_list2_6}, we recall from \cite{Korhonenzhang2023} that $\delta_5(z)\not\equiv0,\pm1,\pm i$. We shall show that $\delta_5(z)$ is a constant. Otherwise, from the proof in \cite{Korhonenzhang2023} we know that $\delta_5(z)$ may be written as $\delta_5(z)=\beta(z)/\alpha_1(z)$, where $\beta(z)$ is a rational function and $\alpha_1(z)$ is in general an algebraic function, namely $\alpha_1(z)$ satisfies an equation of the form $a_2(z)\alpha(z)^2+a_1(z)\alpha(z)+a_0(z)=0$ with some polynomials $a_0(z),a_1(z),a_2(z)$. Moreover, letting $\alpha_2(z)$ also satisfy this equation, we have $\alpha_2(z)/\alpha_1(z)=\delta_5(z)^2$ and $\delta_5(z)$ satisfies an equation of the form $2a_p(z)/\alpha_1(qz)^2=(1+\delta(z))^2/(1+\delta(z)^2)$, where $a_p(z)$ is a nonzero rational function. Together with the relation in \eqref{intro_eq_list2_5 co}, we find that $\beta(z)^2=\alpha_1(z)\alpha_2(z)$ and
\begin{equation}\label{constant assliuoullgui}
\begin{split}
H_2(z)\delta_5(z)^2-H_1(z)\delta_5(z)+H_2(z)=0,
\end{split}
\end{equation}
where $H_1(z)$ and $H_2(z)$ are two relatively prime polynomials such that $H_1(z)/H_2(z)=2a_p(z)^2/\beta(qz)^4$. Then $\delta_5(z)$ is in general an algebraic function defined on two-sheeted Riemann surface. From the relation in \eqref{intro_eq_list2_5 co} we see that $\delta_5(z)$ does not tend to $0$, $\pm i$ and also does not tend to $\infty$ as $z\to\infty$. This implies that the degrees of $H_1(z)$ and $H_2(z)$ are the same. We may write $\delta_5(z)$ as
\begin{equation}\label{constant assliuoullguih0}
\begin{split}
\delta_5(z)=\frac{G_1(z)}{G_2(z)}=\frac{H_1(z)+\sqrt{H_1(z)^2-4H_2(z)^2}}{2H_2(z)}.
\end{split}
\end{equation}
Then we have from the relation in \eqref{intro_eq_list2_5 co} that
\begin{equation}\label{constant ass}
\begin{split}
8\frac{G_1(qz)^4}{G_2(qz)^4}\cdot\frac{G_1(z)^2+G_2(z)^2}{G_2(z)^2}\cdot\frac{G_1(z)}{G_2(z)}=\frac{[G_1(z)+G_2(z)]^4}{G_2(z)^4}.
\end{split}
\end{equation}
Let $z_0$ be a zero of $H_2(z)$ with least order. Then $\delta_5(z)$ has an (algebraic) pole at the point $z_0$. From equation \eqref{constant ass} we see that $qz_0$ is also an (algebraic) pole of $\delta_5(z)$. However, by comparing the order of the poles of the functions on both sides of \eqref{constant ass}, we get a contradiction. Thus $\delta_5(z)$ must be a constant. Then from the discussions in previous section we know that equation \eqref{intro_eq_list2_6} cannot have zero order transcendental meromorphic solutions.

Finally, for equation \eqref{yanagiharaeq13 co}, by the results in \cite{Korhonenzhang2020} we know that the solution $f(z)$ of equation \eqref{yanagiharaeq13 co} is meromorphic and is represented as
$f(z)=\frac{8\lambda(z)^2-(\lambda(z)^2+1)^2}{(\lambda(z)^2+1)^2}$, where $\lambda(z)$ satisfies the $q$-difference Riccati equation in \eqref{riccati  spe2}, which can be transformed to equation \eqref{q-riccati  1}. However, from previous discussion we know that equation \eqref{q-riccati  1} can have only constant solutions when $|q|\not=1$. Therefore, equation \eqref{yanagiharaeq13 co} does not have transcendental meromorphic solutions when $|q|\not=1$. We complete the proof.

\subsection{Growth of the composite function $f(\omega(z))$}\label{growth comparions} 
Let $f(z)$ be a nonconstant meromorphic solution of equation \eqref{first_order_q-ver1} in the autonomous case. In this section we discuss the growth of the composition $g(z)=f(\omega(z))$, where $\omega(z)$ is an entire function such that $\omega(z+1)=q\omega(z)$.

We first look at the equation \eqref{intro_eq_list2_6}. Then the function $g(z)=f(\omega(z))$ satisfies the difference equation
\begin{equation}\label{intro_eq_list2_6 fuq}
g(z+1)^2= \frac{1}{2}\frac{(1+\delta_5)^2}{1+\delta_5^2}\frac{(g(z)-1)(g(z)-\delta_5^2)}{(g(z)-\delta_5)^2},
\end{equation}
where $\delta_5$ is a constant such that $8\delta_5^4(\delta_5^2+1)\delta_5=(\delta_5+1)^4$. Moreover, from the proof in \cite{Korhonenzhang2023} we know that $g(z)$ is written as $g(z)=\text{sn}(\varphi(z))$, where $\text{sn}(\varphi)$ is the Jacobi elliptic function with modulus $1/\delta_5^2$ and $\varphi$ is an entire function of $z$ satisfying the difference linear equation
\begin{equation}\label{intro_eq_list2_6qrwe}
\varphi(z+1) =c_1\varphi(z)+c_2,  \ \ c_1=\pm \frac{1}{2}\frac{(1+\delta_5)^2}{\delta_5^2},
\end{equation}
where $c_2$ is a constant. Then $T(r,\varphi)\geq Lr$ for some $L>0$ and all $r\geq r_0$ with some $r_0\geq 0$.

We estimate the characteristic function $T(r,g)$ using $\text{sn}$ and $\varphi$ as in \cite[Lemma~2.7]{Korhonenzhang2023}. For the composite function $g(z)=\text{sn}(\varphi(z))$, since $\text{sn}(z)$ has positive exponent of convergence of zeros, say $\lambda$, and $\varphi$ has at most one finite Picard's exceptional value, we may choose a constant $w_1>0$ such that $\varphi$ takes in $|z|<t$ every value $w$ in the annulus $w_1<|w|<M(t,\varphi)$, provided that $t$ is large enough. Here $M(t,\varphi)$ denotes the maximum modulus of $\varphi$ on the circle $|z|=t$. Let $\text{sn}$ have $\mu(t)$ zeros in this annulus, counted according to their multiplicity. Then by the definition of $\lambda$, we actually have
\begin{equation}\label{exponent of zeros-00}
\limsup_{r\to\infty}\frac{\log n(r)}{\log r}=\limsup_{t\to\infty}\frac{\log \mu(t)}{\log M(t,\varphi)}=\lambda>0.
\end{equation}
Hence, for some $\tau>0$, there exists a sequence $(t_n)$ such that
\begin{equation}\label{exponent of zeros-01}
\mu(t_n)>\left(M(t_n,\varphi)\right)^{\tau}\geq \left(e^{bt_n}\right)^{\tau},
\end{equation}
where $b$ is a positive constant. The second inequality above follows by the fact that $\log M(t,\varphi)\geq T(t,\varphi)$ for all large $t$. Now, $\text{sn}\circ\varphi$ has at least $\mu(t)$ zeros in $|z|<t$. Making use of \eqref{exponent of zeros-00}, we have
\begin{equation}\label{exponent of zeros02}
\limsup_{r\to\infty}\frac{\log n(r,1/\text{sn}\circ\varphi)}{r}\geq \limsup_{t_n\to\infty}\frac{\log \mu(t_n,1/\text{sn}\circ\varphi)}{t_n}\geq b\tau.
\end{equation}
By the definitions of $N(r,1/g)$ and $n(r,1/g)$ we can deduce that $N(2r,1/g)\geq \frac{1}{2}n(r,1/g)$. Then, by the fact that $T(r,1/\text{sn}\circ\varphi)\geq N(r,1/\text{sn}\circ\varphi)$, we conclude that there is a sequence $(r_n)$ such that the characteristic function $T(r_n,g)$ satisfies
\begin{equation}\label{exponent of zeros0}
\log T(r_n,g)\geq \frac{1}{2}b\tau r_n
\end{equation}
for all large $n$. Similarly, from the discussions in previous section we know that the nonconstant solutions of each of the seven equations \eqref{yanagiharaeq14 co}--\eqref{intro_eq_list2_5} have positive even order or infinite order and thus also gives the estimate in \eqref{exponent of zeros0}.

On the other hand, for the composite function $g(z)=f(\omega(z))$, for any $\varepsilon>0$, we have from \cite[Satz~2.2]{BergweilerJankVol1984} that
\begin{equation}\label{exponent of zeros0uy}
T(r,f(\omega))\leq (1+\varepsilon)T(M(r,\omega)+\omega(0),f),
\end{equation}
for all $r\geq r_1$ and some $r_1>0$. We may choose the entire function $\omega(z)$ to be such that $T(r,\omega)\leq lr$ for some constant $l>0$ and all $r\geq 0$, say $\omega(z)=e^{z\log q}$. Then $\log M(r,\omega)\leq 3T(2r,\omega)$ by an elementary inequality for the relationship between $T(r,\omega)$ and $M(r,\omega)$ (see \cite{Hayman1964Meromorphic}). If $f(z)$ is a zero order transcendental meromorphic solution of equation \eqref{first_order_q-ver1}, then by the definition of order, it follows from the inequality in \eqref{exponent of zeros0uy} that, for any $\varepsilon>0$,
\begin{equation}\label{exponent of zeros0uy1}
\log T(r,g)\leq \log T(M(r,\omega)+\omega(0),f)+O(1)\leq \varepsilon r,
\end{equation}
for all $r\geq r_2$ and some $r_2>0$. This implies that $\log T(r,g)=o(r)$. Comparing this with the estimate in \eqref{exponent of zeros0}, we see that the only possible equations in \eqref{first_order_q-ver1} that can have zero order transcendental meromorphic solutions are \eqref{q-lineareq}, \eqref{q-driccati}, \eqref{yanagiharaeq11 co}, \eqref{yanagiharaeq12 co} or \eqref{yanagiharaeq13 co}.

\section{Concluding remarks}\label{Concluding remarks}

The integrability of difference equations has been demonstrated to be closely related to the existence of meromorphic solutions with minimal hyper-type, as is shown in \cite{AblowitzHalburdHerbst2000,halburdrk:LMS2006,zhengR:18}. The natural integrability criterion for $q$-difference equations is the existence of zero order meromorphic solutions. This is because the transformation $z\to \omega(z)$ with an entire function $\omega(z)$, satisfying $\omega(z+1)=q\omega(z)$, sends $q$-difference equations to difference equations, and the estimate in \eqref{exponent of zeros0uy1} together with the results in \cite{zhengR:18} shows that the existence of sufficiently many zero order meromorphic solutions of a $q$-difference equation implies the proposed difference Painlev\'e property for the related difference equation. However, we should note that the two $q$-difference equations \eqref{yanagiharaeq14 co} and \eqref{yanagiharaeq15 co} do not have any nonconstant zero order meromorphic solutions that are transformed to finite order meromorphic solutions of the two difference equations $f(z+1)^2=(f(z)^2-\kappa)/(f(z)^2-1)$ or $f(z+1)^3=1-f(z)^{-3}$.

Finally, we remark that the transformation $z\to e^{z\log q}$ allows us to consider zero order meromorphic solutions of equation \eqref{first_order_q-ver1} in the punctured complex plane $\mathbb{C}-\{0\}$. Recall the Jacobi elliptic function $\text{sn}(z)=\text{sn}(z,k)$ with the modulus $0<k<1$ and the complete elliptic integral $K$ and $K'$. Then $\text{sn}(z)$ has order $2$ and satisfies the first order differential equation $\text{sn}'(z)^2=(1-\text{sn}(z)^2)(1-k^2\text{sn}(z))$. Let
\begin{equation}\label{Jacobi  1}
\begin{split}
\tau=iK'/K, \ \ \ \ \mathfrak{q}=\exp(\pi i\tau)=\exp(-\pi K'/K).
\end{split}
\end{equation}
Then $\text{sn}(z)$ has the following expression
\begin{equation}\label{Jacobi  2}
\begin{split}
\text{sn}(z)=2\mathfrak{q}^{1/4}k^{-1/2}\sin\frac{\pi z}{2K} \prod_{n=1}^{\infty}\left(\frac{1-2\mathfrak{q}^{2n}\cos \frac{\pi z}{K} +\mathfrak{q}^{4n}}{1-2\mathfrak{q}^{2n-1}\cos \frac{\pi z}{K} +\mathfrak{q}^{4n-2}}\right).
\end{split}
\end{equation}
The above expression can be found in \cite{Baxter1982}. To avoid confusion, in the above expression we have used the notation $\mathfrak{q}$ instead of the usual $q$ there. Look at the $q$-difference equation \eqref{yanagiharaeq14 co}. By the transformation $z\to e^{z\log q}$ we have $f(e^{z\log q})=\text{sn}(z)$ is a solution of the equation $f(z+1)^2=(f(z)^2-\kappa)/(f(z)^2-1)$. Thus
\begin{equation}\label{Jacobi  3 plus}
\begin{split}
f(z)=2\mathfrak{q}^{1/4}k^{-1/2}\sin\frac{\pi\log z}{2K\log q} \prod_{n=1}^{\infty}\left(\frac{1-2\mathfrak{q}^{2n}\cos \frac{\pi\log z}{K\log q} +\mathfrak{q}^{4n}}{1-2\mathfrak{q}^{2n-1}\cos \frac{\pi\log z}{K\log q} +\mathfrak{q}^{4n-2}}\right).
\end{split}
\end{equation}
Recall that $\sin z=\frac{e^{iz}-e^{-iz}}{2i}$ and $\cos z=\frac{e^{iz}+e^{-iz}}{2}$. If the constant $q$ satisfies $\frac{i\pi}{2K\log q} =m$ for some nonzero integer $m$, i.e., $q=\exp(i\pi/2mK)$, then the series in \eqref{Jacobi  3 plus} represents a zero order meromorphic function in the punctured complex plane $\mathbb{C}-\{0\}$ that solves the $q$-difference equation \eqref{yanagiharaeq14 co}.

\vskip  3mm

\noindent{\bf Acknowledgements.} The second author would like to thank professor Rod Halburd of University College London (UCL) for a helpful discussion on the relationship between $q$-difference equations and difference equations during his visit to UCL in~2023.

\end{document}